
\input vanilla.sty
\input definiti.tex
\magnification 1200
\baselineskip 18pt
\overfullrule=0pt
\input mathchar.tex

\define\e{\varepsilon}

\define\pbf{\par\bigpagebreak\flushpar}
\define\a{\alpha}

\def\a{\alpha}

\def\gam{\gamma}

\def\Aut{\hbox{ Aut}}
\def\Om{\Omega}
\def\om{\omega}

\def\varp{\varphi}

\def\vare{\varepsilon}
\def\sig{\sigma}
\def\Sig{\Sigma}

\def\del{\delta}

\def\Im{\hbox{\rm Im}}
\def\Ker{\hbox{\rm Ker \ }}

\def\part{\partial}

\def\cale{{\cal E}}

\def\Area{\text{\rm  Area}}

\title
QUADRATIC EQUATIONS IN GROUPS\\
FROM THE GLOBAL GEOMETRY VIEWPOINT
\endtitle

\author
Alexander Reznikov
\endauthor

\heading{0. INTRODUCTION}\endheading

Consider a group $G$ and a system of equations
$$
\align
&w_1(x_1\ldots x_n)=z_1\\
&\cdots\\
&\vdots\tag *\\
&w_k(x_1,\ldots x_n)=z_k\\
\endalign
$$
where $w_1,\ldots,w_k$ are some words in $x_1,\ldots,x_n$.
It is a classical problem to try to describe the set of solutions of (*) in
$G$.
If the system is quadratic in the sense of [14], then it is essentially
equivalent
to one quadratic equation $w(x_1,\ldots,x_n)=z$, which can be transformed to
one of the canonical forms $$[x_1,x_2]\ldots [x_{2g-1},x_{2g}]=z\tag **$$ or
$$x_1^2\ldots x_g^2=z.\tag ***$$
Even for one homogeneous quadratic equation in a {\bf free} group $G$, the
question of describing all the solutions is nontrivial.
It has been given much effort, starting with the work of Lyndon [17],
see also Lyndon, Wicks [20], Zieschang [39],
Lyndon, McDonough, Newman [19], Culler [6], Chmelevsky [4],
 Burns, Edmunds, Farvouqi [2], Commerford and
Edmunds [5], Goldstein and Turner [10], Rosenberger [27],
[28]. The survey of Schupp [31] describes the state of art by 1980. Recently, a
major
progress has been made,
for $G$ free, in describing the set of solutions in the work of
Piollet [22], Commerford and Edmunds [5],
Grigorchuck, Kurchanov and
Zieschang [10], and Olshansky [23].
The question concerning the solutions of a homogeneous quadratic equation in a
free
group, has been given a satisfactory answer in these papers.
The technique involved uses the cancellation diagrams, described in details
in [18], ch. V. Viewing a free group as a free product of infinite cyclic
groups,
it is quite natural to ask, wether it is possible to study the equations in a
general
free product, provided the information is given about the solutions in each of
the
factors.
The final answer here is due to Rosenberger [27] with previous contribution of
Commerford and Edmonds [5] and Burns, Edmonds and Formanek [3].

In the different direction, Marc Culler [6] has developed a beautiful approach
to
inhomogeneous equations in free groups.
He defined a genus of an element $z\in G'$ similar to the genus norm in the
terminology of Thurston [35] and proved several deep results, e.g.,
on the growth of genus $(z^p)$ as $p\to\infty$.
In the same paper he generalized the finite orbit principle of Chmelevsky [4]
to arbitrary genus, and found nice applications for endomorphisms of free
groups.

Another important contribution has been made by Edmonds [7], see also [8].
This paper contains a decomposition, up to a homotopy, of a continuous map
between
surfaces into a product of a pinch and a branched covering.
This in fact may be transformed to a statement about quadratic equation in
surface groups.
It is also shown in [7], that any degree zero map between closed
surfaces $S,S'$ may be
homotoped to a map from $S$ to a one-skeleton
of $S'$, which in essentially a
homomorphism of $\pi_1(S)$ to a free group.
Thus knowing such homomorphisms, one knows the (homotopy classes of) degree
zero maps
between surfaces.

We also mention the paper of Shapiro and Sonn [32], who have found a very
simple
approach to the Lyndon's rank inequality [18].
They computed the rank of the coupling form in the first cohomology
of a one-relator group, and used the elementary (linear) symplectic geometry
to study epimorphisms on a free group.

Finally, a very general picture of algorithmic aspects for equations in
free groups has
been given in Razborov [24](and, with more detailed information, in his
thesis).

The present paper suggests somewhat unexpected approach to equations in groups
via
deeper parts of Riemannian geometry of manifolds.
The central analytic tool that we will use is an existence theorem for
minimal surfaces in Riemannian manifolds
satisfying some uniformity conditions generalizing well-known theorem of
Sacks-Uhlenbeck [29] and Schoen-Yau [30].
It uses classical genus reduction process as described by Ville [36].
In our previous paper [25] we have used this kind of argument
for proving Thurston-Gromov finiteness theorem for
(geometrically) hyperbolic groups.

The power of the approach presented here allows us to cover nearly all the
results,
described above, on quadratic equations in free groups,
homogeneous and inhomogeneous and to give far-going generalizations to free
products and subgroups of geometrically hyperbolic groups.
There will be given two different proves of the Classification Theorem for
solutions
of homogeneous quadratic equations in a free group.
One fits a more general framework of the Rosenberger's Classification Theorem
for
equations in free products, which we approach from our point of view in the
section 4.

The other proof is based on the work of Schoen and Yau concerning the minimal
surfaces
in compact three-manifolds of positive scalar curvature; it will be presented
in the
section 3.
It applies to a very restricted class of free products (including all free
groups)
but gives much more precise information on the structure of solutions;
it turns out that all solutions are ``elementary'' in a sense explained below
in 3.2.

We then state and prove the Classification Theorem for solutions of quadratic
homogeneous
equations in geometrically hyperbolic groups.
These are fundamental groups of complete manifolds of pinched negative
curvature, satisfying the uniformity
and growth conditions of 5.1.
This is a subclass in the
class of Gromov word-hyperbolic groups.

Our existence theorem on minimal surfaces reduces, in the case when the
ambient space is a surface, to a version of the  well-known result of
Edmunds on decomposition of maps between surfaces, given in a very
precise and geometric form (our branched convergings are actually conformal
maps).

Next we turn to the inhomogeneous equations.
This corresponds to the solutions of Pleateau problem in the minimal surfaces
theory.
Using the current (de Rham) norm in the loop space we recover the result
of Culler on genus $(z^p)$ in much more general context of subgroups of
hyperbolic
groups, as far as his generalizations of Chmelevsky' finite orbit principle.
To apply our technique to actually free groups, we need to realize these as
commutators subgroups of excellent fiber knot groups and to use in part the
Thurston hyperbolization theorem for knot manifolds [34].

This work owes its existence to the help of many people.
It was started in the most warmful and stimulating atmosphere during my visit
to Ruhr-Universit\"at Bochum in 1991--1992.
The warmest thanks are due to Heiner Zieschang, Martin Lustig, Frank Levin,
J\"urgen Jost, Shicheng Wang as well as to Ilya Rips and Marina Ville for
valuable discussions.


\heading{1. PRELIMINARIES ON RIEMANNIAN SURFACES}\endheading

1.1. Most of the essential results of this paper equally holds
for nonorientable surfaces as well as for orientable.
For the reader's convenience we will adopt the following strategy:
complete proof will be provided for the orientable case with the precise
specification
of the changes necessary to cover the case of nonorientable surfaces.
The formulation of the main classification results are more complicated for the
nonorientable case, basically because of the more involved ierarchy of simple
closed curves on a nonorientable surface, see [38].

We begin by the description of the genus reduction process of Douglas-Courant.
The presentation given here follows that of [36], [26].

\proclaim{1.2. Lemma} Let $\Sig^g$, $g\ge 2$ be an orientable surface and let
$X$ be
a CW complex.
Let $f:\Sig^g\to X$ be a continuous map.
Then one of the following two possibilities holds:
\item{\rm (i)} $f_*:\pi_1(\Sig^g)\to \pi_1(X)$ is essentially injective

\item{\rm (ii)} There exists a homotopy decomposition
$$
\matrix
\Sig^g & \mathop{\longrightarrow}\limits^f& X\\
\searrow && \nearrow\\
&\Sig^{g_1}\bigvee \Sig^{g_2}&\\
\endmatrix
$$
with $2-g=(2-g_1)+(2-g_2)+2$, and $f_*$ decomposes as
$$
\matrix
\hat \pi_1(\Sig^{g_1})_{\underset \bbz\to *}\hat\pi_1(\Sig^{g_2})=\pi_1
(\Sig^g)&\mathop{\longrightarrow}\limits^{f_*}
&\pi_1(X)\\
\qquad\qquad\qquad\qquad\qquad\qquad
\searrow&&\nearrow\\
&\pi_1(\Sig^{g_1})*\pi_1(\Sig^{g_2})&\\
\endmatrix
$$
where $\hat \pi_1(\Sig^{g_i})$ is the (free) fundamental group of $\Sig^{g_i}$
punctured at one point.
\endproclaim

\demo{Proof} Recall that   $f_*$ is not essentially injective, means that
$\Ker f_*$
contains a class of a simple closed curve $\gam$.
If this curve $\gam$ is separating, we may pinch it to a point
as shown in Fig 1.3 (left).
The resulted space is a wedge and (ii) follows. If
$\gam$ is not separating then([38])
it is a meridian of a handle as in the Fig 1.3 (right).
Let $\del$ be the longitude of the handle and $\vare$ be the attaching circle.
Then, up to conjugacy in $\pi_1(\Sig^g)$, $\vare=[\gam,\del]$ so
$f_*\vare=[1,f_*\del]=1$.
Since $\vare$ is separating, we may proceed as before.
\demo{1.4. Definition} Let $X$ be a CW-complex.
let $H:\pi_2(X)\to H_2(X)$ be the Hurewitz map, and
let
$\overline{\pi_2(X)}=\Im H$.
Let $z\in H_2(X)/\overline{\pi_2(X)}$.
The Thurston genus norm of $z$, denoted $||z||_g$ is the minimal absolute
value of the Euler characteristic of a singular surface
$f:\Sig^g\to X$, representing $z$.
It is allowed for the orientable surface $\Sig^g$ to be disconnected.

\proclaim{1.5. Lemma} Let $X$ and $z$ be as above.
Suppose $f:\Sig^g\to X$ is a continuous map such that
$\chi(\Sig^g)=||z||_g$.
Then the restriction of $f$ on every connected component of $\Sig^g$ is
essentially
injective.
\endproclaim

\demo{Proof} See [26], lemma 1.

The following standard fact will be frequently used in the paper,
c.f.Culler
[6].

\proclaim{1.6. Lemma} Let $X$ be a CW-complex and let
$[\gam]\in [\pi_1(X),\pi_1(X)]$.
Then there exists a map $f$ of a  surface $\Sig^g$ with a removed disc
$D^2\subset \Sig^g$  to $X$, such that $f|_{\part D_2}=\gam$.
\endproclaim

\demo{Proof} Let $[\gam]=a_1b_1a_1^{-1}b_1^{-1}\ldots a_mb_ma_m^{-1}b_m^{-1}$
in
$\pi_1(X)$.
The fundamental group of $\Sig^m\setminus D^2$ is free  on generators
$x_1,y_1,\ldots x_m,y_m$, and $[\part D^2]=[x_1,y_1]\ldots [x_m,y_m]$.
Any map $f:\Sig^m\setminus D^2$ such that $f_*(x_i)=a_i$ and
$f_*(y_i)=b_i$ may be homotoped to yield $f|_{\part D^2}=\gam$.

Now we specify the changes in lemma 1.3., to yield the nonorientable case.

\proclaim{1.7. Lemma} Let $\Sig^g$ be a nonorientable surface and let
$f:\Sig^g\to X$ be a continuous map.

Then one of the following possibilities hold:
\item{\rm (i)} $f_*$ is essentially injective

\item{\rm (ii)} $f$ may be homotoped to a decomposition
$$
\matrix
\Sig^g & \mathop{\longrightarrow}\limits^f& X\\
\searrow && \nearrow\\
&\Sig^{g_1}\bigvee \Sig^{g_2}\\
\endmatrix$$
with $\chi(\Sig^g)=\chi(\Sig^{g_1})+\chi(\Sig^{g_2})-2$, and $f_*$ decomposes
as
$$
\matrix
\pi_1(\Sig^g)=\hat \pi_1(\Sig^{g-1})_{\buildrel * \over \bbz}
\hat \pi_1(\Sig^{g_2})& \longrightarrow& \pi_1(X)\\
\hbox{\hskip6truecm}\searrow &&\uparrow\\
&& \pi_1(\Sig^{g_1})*\pi_1(\Sig^{g_2})\\
\endmatrix
$$
where $\hat \pi_1(\Sig^{g_i})$ is the fundamental group of $\Sig^{g_i}$
punctured at one point
\item{\rm (iii)} There exists a pair $(\Sig^{g'},f')$ consisting of an
orientable
surface $\Sig^{g'}$, and a map $f':\Sig^{g'}\to X$ such that
$(\Sig^g,f)$ is (up to homotopy) obtained from $(\Sig^{g'},f')$ by the
following
surgery:
remove a small disc $D^2\subset \Sig^{g'}$ such that $f'|_{D^2}=\text{\rm
const}$,
glue the opposite points of $\part D^2$ together to obtain $\Sig^g$,
and define $f$ by restriction.

\item{\rm (iv)} There exists a pair $(\Sig^{g'},f')$ consisting
a surface $\Sig^{g'}$ and a map $f':\Sig^{g'}\to X$ such
that $(\Sig^g,f)$ is obtained from
$(\Sig^{g'},f')$ by the following surgery: remove two discs
$D_1,D_2\subset \Sig^{g'}$ such that $f'|_{D_i}=\text{\rm const}$, glue the
boundaries of $D_1$ and $D_2$ together to obtain $\Sig^g$, and define $f$ by
restriction.

\endproclaim

\demo{Proof} Suppose $f$ is not essentially injective and
$[\gam]\in\Ker f_*$ is a simple closed curve.
If $\gam$ is separating, then (ii) holds.
If $\gam$ is not separating, then either $2\gam$ is
homotopic to a separating curve [38], (3.5.7) and
we proceed as before with $\gam$ replaced by $2\gam$, of one of
the possibilities, described in [38], (3.5.8)--(3.5.10) holds.
This corresponds to the cases (iii) and (iv) o
f the lemma.



\heading{2. THE FUNDAMENTAL EXISTENCE THEOREM FOR MINIMAL SURFACES}\endheading

2.1. For a Riemannian manifold, $M$, we adopt the {\bf uniformity and
the growth conditions} in the form: the injectivity radius
$i_x(M)\to \infty$ as $x\in M$ escapes
all compact subsets of $M$ and for some $R>0$ and for any $x\in M$
there exists a uniformly by $x$ bi-Lipschitz diffeomorphism of the
Euclidean ball $B(0,R)$ on a neighborhood of $x$.
The Morrey regularity theory implies, as in [30],  that any $W^1_2$ energy
minimizing map of
a compact surface to $M$ will be smooth, hence harmonic.

Let $\Sig^g$, $g\ge1$, be an orientable closed surface, and let $\gam$ be
a simple closed curve in $\Sig^g$.

\demo{2.2. Definition} An elementary pinch is a continuous map $\del$,
defined as follows.
If $g\ge 2$, then we assume that $\gam$ is separating, and $\del$ is the map
$$
\del:\Sig^g\to\Sig^{g_1}\bigvee\Sig^{g_2}
$$
of $\Sig^g$ on the wedge of two closed orientable surfaces, obtained
by contracting $\gam$ to a point.
If $g=1$, then
$$
\del:\Sig^1\to S^1
$$
is defined by contracting $\gam$, viewed as a meridian of the torus
$\Sig^1$ to a point followed by a retraction on a longitude.

\demo{2.3 Definition} A pinch
$\del:\Sig^g\to\Sig^{g_1}\bigvee\cdots\bigvee \Sig^{g_r}\bigvee S^1_{(1)}
\bigvee \cdots S^1_{(s)}$ is defined inductively as composition of elementary
pinches.

\proclaim{2.4. Theorem} Let $\Sig^g$ be a closed oriented surface,
and let $M$ be a Riemannian manifold which satisfies the uniformity conditions
and the growth conditions.Suppose that $M$ is either compact or of pinched
non-positive curvature.
Let $f:\Sig^g\to M$ be a continuous map.
Then there exists a pinch $\del:\Sig^g\to \Sig^{g_1}\bigvee\cdots\bigvee
\Sig^{g}_r\bigvee S^1_{(r+1)}\bigvee\cdots S^1_{(r+s)}$,
a collection of (branched) minimal immersions $\varp_i:\Sig^{g_i}\to M$,
$1\le i\le r$, and a collection of continuous maps
$\hat \varp_j:S^1_{(j)}\to M$, $1\le j\le r$, such that the diagram
$$
\matrix
\pi_1(\Sig^g)&\mathop{\longrightarrow}\limits^{f_*}&\pi_1(M)\hbox{\hskip3truecm}\\
\\
\searrow & &
\uparrow \hat\varp_{1_*}*\hat\varp_{2_*}*\cdots*\hat\varp_{(r+s)_*}\\
\\
&\pi_1(\Sig^{g_1})*\cdots *\pi_1(\Sig^{g_r})*\bbz*\cdots *\bbz\\
\endmatrix
$$

is commutative.
Here $\hat\varp_i$, $1\le i\le r$ are some homomorphisms, conjugate to
$\varp_i$ in $\pi_1(M)$.
If $\pi_2(M)=0$, then $f$ is freely homotopic to $\bigvee \varp_i$.
\endproclaim

\demo{Proof} We begin by applying the lemma 1.3. possibly several times,
to find a decomposition 2.4  where all $\varp_i$, $1\le i\le r$
are essentially injective maps.
We need to show that
 $\varp_i$ may be replaced by minimal
immersions $\tilde\varp_i$ which induce
up to conjugation in $\pi_1(M)$ the same homomorphisms of the
fundamental groups.
The lines which follow actually mimic the argument of Sacks-Uhlenbeck [26]
and Schoen-Yau [30].
Fix metrics $h_i$ on $\Sig^{g_i}$. Using the direct method of variational
calculus, we
find $W^1_2$ energy minimizing maps $\psi_i:\Sig^{g_i}\to M$ inducing (up to
conjugacy)
the same map of fundamental groups, as $\varp_i$.

By the regularity theory of Morrey, all $\psi_i$ are smooth, hence harmonic.
Next, we let the conformal classes $[h_i]$ of the metrics $h_i$ vary inside the
Teichmuller space
$T_{6g_i-6}$.
Let $[h_i^{(\a)}]$ be the sequence of conformal classes,
such the $Energy \ (\psi_i^{(\a)})$ decreases to its infinum taken over
$T_{6g_i-g}$.
Because of the conformal invariance of the energy, we may assume
$h_i^{(\a)}$ to be the (uniquely defined in their conformal classes)
hyperbolic metrics.

Then by the Collar theorem of  we conclude, as in [27],
that the minimal length of a closed geodesic on $(\Sig^{g_i},h_i^{(\a)})$
is bounded away from zero.
Observe that here we need the uniformity conditions on $M$, to assume that
the lengths of essential curves of $M$ are bounded away from zero, and the
essential
injectivity of $(\varp_i)_*$.
By the Mumford compactness criterion, the sequence $[h_i^{(\a)}]$
contains a subsequence, which we immediately relabel
$[h_i^{(\a)}]$, whose image in the moduli space $M_{6g-6}$
converges.
Twisting by a diffeomorphism $\varp^{(\a)}_i:\Sig^{g_i}\to \Sig^{g_i}$
we arrive to a converging sequence, again relabeled $[h_i^{(\a)}]$, already in
the Teichmuller space itself.Now, the estimates of [29] show that
$\psi^{(\alpha)}_i$ are H\"older equicontinious. Moreover,
$Im(\psi^{(\alpha)}_i)$ cannot escape all compact sets in $M$, since otherwise
$\psi^{(\alpha)}_i$ could be lifted to the universal cover of $M$, which is
imposible by the maximum principle [ ]. So $\psi^{(\alpha)}_i$ contains a
$C^0$-converging subsequence,whose limit is minimal by [29],[30].

If $\pi_2(M)=0$ then $f$ is freely homotopic to the wedge $\bigvee \varp_i$.
In particular, if $M$ is a closed surface, we recover the following well-known
results of Edmonds [7] (the Main Theorem)

\proclaim{2.5. Theorem (Edmonds)} Any continuous map $f$ between surfaces
$\Sig$ and $M$ can be decomposed, up to a homotopy, as a pinch followed by
a branched covering.
\endproclaim

Another consequence is a theorem of M. Ville [36]:
any element in $H_2(M)/\overline{\pi_2(M)}$
can be represented by a minimal surface.
Here $\overline{\pi_2(M)}$ is the image of $\pi_2(M)$ in $H_2(M)$ under the
Hurewitz
homomorphism.
Combining this with the theorem of Sacks-Uhlenbeck,
we get that any element in $H_2(M)$ is represented by a minimal surface.
We refer to our paper [26] for various geometric applications of this
result.

Observe that the definition of a pinch, given in [7], differs slightly from
ours.


\heading{THE THREE-DIMENSIONAL SURGERY AND HOMOGENOUS QUADRATIC}\endheading
\heading{EQUATIONS IN FREE PRODUCTS OF SPECIAL TYPE}\endheading

3.1. In this section, we apply the Existence Theorem for minimal surfaces to
the
study of the solutions of the quadratic equation (*) is a free product of
groups specified
below in 3.3.
The idea is to realize such a group as a fundamental group of a compact
three-manifold of positive scalar curvature and then to apply the Existence
Theorem 2.4.
On the other hand, the result of Schoen-Yau [30] states that there are
no stable minimal surfaces of genus more then zero in such
three-manifold.
This will mean that the pinch $\del$ in 2.4 acts actually to a wedge of
circles, which corresponds to {\bf elementary}
homomorphism of the surface group, in the sense of the following
definition.
Throughout this section $\Sig^g$ stands for a closed {\bf orientable}
surface.

\demo{3.2. Definition} Let $G$ be a group and let
$\rho:\pi_1(\Sig^g)\to G$ be a homomorphism.
We will call $\rho$ to be an {\bf elementary homomorphism} if there
exists a set of canonical generators $(u_i,v_i)_{i=1}^g$ of
$\pi_1(\Sig^g)$ such that $\rho(u_i)=1$ (then $\rho(v_i)$ may be completely
arbitrary elements in $G$).

\proclaim{3.3. Theorem} Suppose $G=\mathop\ast\limits_{i=1}^r G_i$,
where $G_i$ runs through the following list:
$\bbz_m$, $\bbz_m$, $m\ge 2$, $\bbz_m\times D^*_n$, $(m,n)=1$, $\bbz_m\times
T^*_n$, $(m,6)=1$, $\bbz_m\times S_4$, $(m,6)=1$, $\bbz_m\times A_5$,
$(m,30)=1$.
Then any homomorphism $\rho:\pi_1(\Sig)\to G$ is elementary.

Equivalently, any set of solutions of (**) is obtained from an elementary
solution
($x_i=1, y_i$ arbitrary) by a repetition of the Dehn-Nielsen transformations.
\endproclaim
Here $T_n$ stands for the group $\{ x,y,z|x^{3^n}=y^4=1, y^2=z^2, xyx^{-1}=y,
xzx^{-1}=yz, yzy^{-1}=z^{-1}$
The theorem 3.3. contains, as a special case $G=F_r$, the main result
of Zieschang [39], Piollet [22], and Grigorchuk-Kurchanov-Zieschang [10].

\proclaim{3.4. Lemma} Let $G$ be a group from the list of the
theorem 3.3.
Then $G$ can be realized as a fundamental group of a three-manifold of positive
scalar curvature.
\endproclaim

\demo{Proof} Observe that $\bbz=\pi_1(S^2\times S^1)$ and $S^2\times S^1$
can be given a product metric of $(S^2,can)\times (S^2,can)$ with positive
scalar
curvature.

Next, all other groups $G_i$ act freely on $S^3$ [37], hence $S^3/G_i$ admits
desired metric.
To deal with free products, we recall a strong result of Gromov-Lawson
and Schoen-Yau, which says that positive scalar curvature manifolds admit
surgery in
codimension at least three.
In particular, a connected sum of three-manifolds of positive scalar
curvature again carries a metric of positive scalar curvature.
This proves the lemma, since the fundamental group of a connected sum is a free
product.

\demo{Proof of the theorem 3.3} Let $M$ be a three-manifold with fundamental
group $G$, which exists according to the previous lemma.
Let $f:\Sig^g\to M$ be a map which induces the homomorphism
$\varp$.
By the Theorem 2.4 we know that $\rho$ decomposes as $*\varp_{i_*}\circ\del_*$,
where $\varp_i$, $1\le i\le r$, are stable minimal maps
of surfaces of genus $\ge 1$, to $M$.
However, by a result of Schoen-Yau [30], theorem a three-manifold of positive
scalar
curvature carries no stable minimal surface.
Hence $r=0$ and $\rho$ decomposes as a pinch
$\vee\to \vee S^1_{(i)}$ followed
by a map of $\vee S^{(1)}_{(i)}$ to $M$.
This precisely means that the homomorphism $\rho$ is elementary.

\demo{3.5. Remark} If $G$ is a fixed finite group, then, of course,
the number of homomorphisms $\rho:\pi_1(\Sig^g)\to G$ is finite for any
fixed $g$.
However as $g$ grows, it becomes a very
difficult problem to describe all these
homomorphisms (it is much easier to find their number, [21]).
We note that if $H_2(G,\bbz)\not=0$, then there are homomorphisms which
are not elementary in the sense of 3.2.
The converse to this statement is wrong.

\demo{3.6. Remark} Recall that a solution $(x_i,y_i)$ of (**) in $G$
is called free, if the corresponding homomorphism $f:\pi_1(\Sig)\to G$
factors through a free group [18].
Any elementary solution is clearly free.
The converse is also true and follows from 3.3.

So the statement of 3.3 may be reformulated in the following way:
for $G$ as in 3.3. all solutions of (**) are free.


\heading{4. THE CLASSIFICATION THEOREM FOR SOLUTIONS OF HOMOGENEOUS}\endheading
\heading{QUADRATIC EQUATIONS IN FREE PRODUCTS}\endheading

4.1. Here we address the problem of an  algebraic description for
solutions of homogeneous quadratic equations (**), (***)( in a general free
products $G_1*G_2$ in terms of solutions in each of $G_i$, $i=1,2$.
See [] [] for some partial result for small genus.
We will prove here the general results of Rosenberger [27]  below,
which resolves this problem completely.
For simplicity we work with finitely presented groups $G_i$.
The obvious induction argument extends our result to an arbitrary
number of groups.
The idea of our approach is to realize $G_i$ as a fundamental group of a
compact $5$-manifold $M_i$ and then to study maps of surfaces
$\Sig^g$ to $M_1\#M_2$ using the tranversality.
We have a good reason to believe that this approach will prove very
useful in studying similar problems in relations to amalgamated
products, where the connected sums will be replaced by a surgery along a
submanifold.

\proclaim{4.2. Theorem} Let $\rho$ be a homomorphism from a surface group
$\pi_1(\Sig^g)$, $\Sig^g$ oriented, to $G$.
There exists a collection of homomorphisms
$\rho_i:\pi_1(\Sig^{g_i})\to G_1$,  $1\le i\le r_1$,
$\rho_i:\pi_1(\Sig^{g_i})\to G_2$,
$r_1+1\le i\le r_1+r_2$, and $\rho_i:\bbz\to G$,
$r_1+r_2+1\le i \le r_1+r_2+s$, a pinch
$\del:\Sig^g\to \vee\Sig^{g_i}\vee S_{(i)}^1$,
and a collection of elements $z_i$, $1\le i\le r_1+r_2+s$ in $G$,
such that $\rho$
decomposes as
$$
\matrix
\pi_1(\Sig^g)&\mathop{\longrightarrow}\limits^\rho& G\\
\hbox{\hskip2truecm}\searrow \del_*&&\nearrow *(z_i\rho_iz_i^{-1})\\
&\pi_1(\Sig^{g_1})*\cdots*\pi_1(\Sig^{g_{r_1+r_2}})*\bbz*\cdots*\bbz\\
\endmatrix
$$
\endproclaim

\proclaim{4.3. Lemma} Let $G$ be a finitely presented group.
There exists a compact $5$-manifold, $M$, such the $\pi_1(M)=G$ and
$\pi_2(M)=0$.
\endproclaim

\demo{Proof of the theorem 4.2}
We begin by realizing $G_i$ as a fundamental group of a compact
$5$-manifold $M_i$ with $\pi_2(M_i)=0$, as above.
Let $M=M_1\#M_2$, then $G=\pi_1(M)$.
Given $\rho:\pi_1(\Sig^g)\to G$, fix $f:\Sig^g\to M$ with $f_*=\rho$,
up to a conjugacy in $G$.
We denote by $S$ a separating $4$-sphere in $M$.
First we note that $\pi_2(M)=0$.
This actually will follow from the argument below.
Assuming this, we use 1.3. to deform $f$ to a composition of a pinch $\rho$
followed by a wedge of maps $\varp_i:\Sig^{g_i}\to M$, $1\le i\le r$ and
$\varp_i:S^1_{(i)}\to M$, $r+1\le i\le r+s$, such that $\varp_i$, $1\le i\le r$
are essentially invective.
We claim that any $\varp_i$, $1\le i\le r$ is freely homotopic to a map,
whose image lies in one of $M_i$.
Indeed, smooth $\varp_i$ and make it transversal to $S$.
This is possible since $S$ is of codimension one in $M$.
Let $\mu(\varp_i)$ be the number of connected components of
$\varp_i^{-1}(S)$.

We assume $\mu(\varp_i)$ is minimal possible in the
homotopy class of $\varp_i$.
If $\mu(\varp_i)=0$, then the image of $\varp_i$ lies in one of $M_i$,
and we are done.
So we may assume $\varp_i^{-1}(S)$ be nonempty.
Let $C\subseteq\varp_i^{-1}(S)$ be a connected component of $\varp_i^{-1}(S)$.
If $[C]=0$ in $\pi_1(\Sig^{g_i})$ then $C$ bounds a disc, $D$, and we may
assume
that no other components of $\varp_i^{-1}(S)$ are in $D$, otherwise we will
replace $C$ by one of these components.
Then $\varp_i$ maps $D$ to one of $M_i$, say $M_1$,
and we can push $\varp_i|_D$ out of $M_1$ (here we need $\pi_2(M_i)=0$),
reducing $\mu(\varp_i)$.
This argument also shows that $\pi_2(M)=0$.
So $C$ is an essential closed curve,
and $\varp_i$ is not essentially injective, a contradiction.
So $\varp_i$ may be homotoped to a map, whose image is in one of $M_i$.
Back to the homomorphisms, this gives precisely the decomposition of
4.2. Q.E.D.

As an immediate corollary, we get another proof to the Classification Theorem
for equations in free groups, c.f. 3.3.

4.5. We now turn to the study of homomorphisms of a nonorientable
surface groups to $G_1*G_2$, which corresponds to the equation (***).
The crucial point in the proof of the Classification Theorem 4.2. below
is contained in the following lemma.

\proclaim{Lemma} Suppose $f:\Sig^g\to M_1\#M_2$, is a map of an nonorientable
surface
of the genus $g$. Suppose $f$ is not freely homotopic to a map,
whose image is contained in one of $M_i$.

Then $\Ker f_*$ contains a class of an essential closed separating curve.
\endproclaim

\demo{Proof} First we deform $f$ to a smooth map which is transversal to $S$,
and relabel the new map by $f$ again.
The preimage $f^{-1}(S)$ is nonempty and is a disjoint union of simple closed
two-side curves, say $C_1,\ldots ,C_m$.
As in the proof of 4.2., we may think of all these curves to be
essential.
The homology class $[C_1]+\cdots+[C
_m]$ is zero in
$H_1(\Sig,\bbz_2)$ hence if $m=1$, then the (unique) curve $C_1$ is
separating, and we are done.
So we assume that $m\ge 2$ and non of $C_i$ is separating.
Cut $\Sig$ along $C_1$.
There are two possible cases:
\item{(i)} The resulting surface is nonorientable.
Then by [38], $C_1$ has representation $V_1V_2$ in $\pi_1(\Sig^g)$, where
$V_1,\ldots ,V_g$ is a set of generators for $\pi_1(\Sig^g)$ with the canonical
relation
$V_1^2\cdots V_g^2=1$.

Since $f(C_i)$ is contained in the sphere $S$, we have obviously
$f_*(C_i)=1$ in $G$, so $f(V_1)f(V_2)=1$, hence $f(V^2_1V_2^2)=1$.
The element $V_1^2V_2^2$ is represented by a closed separating curve in
$\Sig^g$ and is contained in $\Ker f_*$, as stated.

\item{(ii)} The resulting surface $\tilde\Sig$ is orientable.
Then $\Sig^g$ is obtained by a surgery, described in 1.7(iv).
Consider $C_2$ as a curve in $\tilde\Sig$.
If it is not separating, then it is a meridian of a handle,
and we immediately get a closed separating curve (and attaching circle of the
handle) representing an element of $\Ker f_*$, as in 1.3.
So we may assume $C_2$ is separating (Fig 4.5).
In the case of Fig 4.5,right, the curve $C_2$ is separating in $\Sig$.
So it is enough to study the case of Fig 4.5, left.

The surface $\Sig^g$ is therefore a Klein bottle
with several handles and M\"obious
bands attached.
We may localize these handles and M\"obious bands
within a disc $D^2$ in the Klein bottle.
We claim the boundary circle $\part D^2$ is a curve we need.
Indeed, if is a closed separating essential curve, and we need
only to prove that $[\part D^2]\in \Ker f_*$. Let $C_3$ be the longitude
of the Klein bottle, then
$[\part D^2]=[C_3][C_2]C_3^{-1}][C_2]$, and since $f_*[C_2]=1$, we get
$f_*[\part D^2]=1$. Q.E.D.

4.6. We may now switch on the decomposition procedure of 1.7 using the previous
lemma, and to produce the canonical form of the homomorphism $f_*$ is the
nonorientable case.

\demo{4.7. Definition} Let $\Sig$ be a surface (possibly nonorientable).
An elementary pinch is either
a map, defined by 2.2, if $\Sig$ is not a Klein bottle,
or a projective plane, or the (homotopically unique) map $f:\Sig^2\to \Sig^1$
of the Klein bottle to $S^1$, which induces the generator of
$H^1(\Sig^2,\bbz)=\bbz$,
or the map $\bbr P^2\to pt$, or the map $\Sig^2=\bbr P^2\#\bbr P^2
\mathop{\longrightarrow}\limits_{id\ \vee\text{ const}}\bbr P^2$.

A pinch $f:\Sig\to\Sig^1\vee\cdots\vee\Sig^k\vee S^1\cdots \vee S^1$ is a
composition of elementary pinches.

4.8. Now we state the Classification Thoerem for homomorphisms
$\rho:\pi_1(\Sig)\to G_1*G_2$.

\proclaim{Theorem} Any homomorphism $\rho:\pi_1(\Sig)\to G_1*G_2$ admits a
decomposition
$$
\matrix
\pi_1(\Sig) & \longrightarrow
 & G_1*G_2\\
\searrow\del_*&&\nearrow *(z_ip_iz_i^{-1})\\
&\pi_1(\Sig^{g_1})*\cdots *\pi_1(\Sig^{g_{r_1+r_2}})*\bbz*\cdots*\bbz\\
\endmatrix
$$
where $\del$ is a pinch in the sense of 4.7., and $\rho_i$, $z_i$ are as in
4.2.
\endproclaim

\demo{Proof} The only case which is not yet covered, concerns maps of the Klein
bottle
$\Sig^2$ to $M$ which are not essentially injective.
Suppose $f:\Sig^2\to M$ is such a map and $\gam$ is an essential simple
closed curve such that $[\gam]\in\Ker f_*$.
 Then by the classification of [38], we have the following cases.
\item{(i)} $[\gam]=V_1$. Then homotopically $f$ factors through the pinch
$\Sig^2\to\bbr P^2$.

\item{(ii)} $[\gam]=V_1V_2$. Then homotopically $f$ factors through the pinch
$\Sig^2\to S^1$.

\item{(iii)} $[\gam]=V_1^2$. Then homotopically $
f$ factors through the pinch
$\Sig^2\to\bbr P^2\bigvee\bbr P^2$.

\heading{5. CANONICAL STRATIFICATIONS IN MODULI SPACES OF
SOLUTIONS}\endheading
\heading{TO HOMOGENEOUS QUADRATIC EQUATIONS IN}\endheading
\heading{GEOMETRICALLY
HYPERBOLIC GROUPS}\endheading

5.1. Throughout this section $G$ is a fundamental group of a
Riemannian manifold $M$ with a sectional curvature $K(M)$ pinched as
$-K\le K(M)\le -k<0$ and satisfying the uniformity conditions
and the growth conditions (e.g., compact).
Such groups will be called geometrically hyperbolic.
Observe that a f.g. free group is not a fundamental group of
a compact manifold of negative curvature, but is a fundamental group of
a manifold, satisfying the conditions above, e.g., an infinite cyclic
covering of an excellent knot complement, see 5.6.

It turns out  that the Classification Theorem for quadratic equations in
free groups in a form proved in 3.3 an 4.2 is ``almost true'' for
geometrically hyperbolic groups.
To make it precise, we need a following definition.

\demo{5.2. Definition} Let $\rho:\pi_1(\Sig^g)\to G$, $\Sig^g$ oriented,
be a homomorphism.
We say that $\rho$ has a defect at least $d$, denoted
$def(\rho)\ge d$, if $\rho$ admits a decomposition 2.4  with $s\ge d$.
A maximal $d$ such that $def(\rho)\ge d$ is called the
defect of $\rho$.

\proclaim{5.3. Theorem} Let $G$ be a geometrically hyperbolic group.
Consider all homomorphisms $\rho:\pi_1(\Sig^g)\to G$, of defect $d$,
so that $\rho$ admits a decomposition
$$
\matrix
\pi_1(\Sig^g)
&\longrightarrow&G\qquad\qquad\\
\hbox{\hskip2truecm}\del_*\searrow&&\uparrow\varp_1*\cdots*\varp_{r+d}\\
\\
&&\pi_1(\Sig^{g_1})*\cdots *\pi_1(\Sig^{g_r})*\bbz*\cdots*\bbz\\
\endmatrix
$$
with $\del$ a pinch and $\varp_i$, $1\le i\le r$, essentially injective.
Then for a fixed $\del_*$, there are only finite number of possibilities
up to an automorphism of $\pi_1(\Sig^{g_i})$
and conjugation in $G$ for $\varp_i$, $1\le i\le r$.
\endproclaim

\demo{Proof} Granted 2.4
, it is enough to show that for a fixed $g$,
there are only finite number of conjugate classes of essentially
injective homomorphisms $\varp:\pi_1(\Sig^g)\to G$.
This is a generalization of the well-known theorem of Thurston [35] and
Gromov [13], where the actual
{\bf injectivity} was demanded.
In our previous paper [25], we gave an analytic proof, which in fact applies
here
without any changes.
Indeed, consider a sequence of nonconjugate homomorphisms
$\varp_i:\pi_1(\Sig^g)\to\pi_1(M)=G$, and let $f_i:\Sig^g\to M$
be maps which induces $
\varp_i$ up to conjugacy.
Fix a metric $h$ of curvature $-1$ on $\Sig^g$ and consider the harmonic
maps $\bar f_i:\Sig^g\to M$, inducing $\varp_i$.
This is possible by Sacks-Uhlenbeck [29] and Schoen-Yau [30].
Let $E_i^h(\bar f_i)$ be the energy of $\bar f_i$.
By the Thurston-type inequality of [25], $Area\ (\bar f_i)$ stays bounded as
$i\to\infty$,
so we may
perturb $h$ to some $h_i$, again hyperbolic, such that $E_i^{h_i}(\bar f_i)$
is bounded.
Now the argument of [30] shows that the class of $h_i$ in the moduli space
$M_{6g-g}$ stays in some compact set (here we need the uniformly conditions
for $M$), so we can replace $h_i$ by some fixed $\bar h$, twisting
$\bar f_i$ if necessary by a diffeomorphism of
$\Sig^g$ such that $E_i^h(\bar f_i)$ remains bounded.
Applying the collar argument of [30] again, we see that for a given
curve $C$ in $\Sig^g$, the minimal length of curves in the homotopy class of
$\bar f_i(C)$ stays bounded.
Since $M$ satisfies growth conditions, there are but finite number of conjugacy
classes in $G$ which may be images of $[C]$ under $f_{i_*}$.
Since $\pi_1(\Sig^g)$ is not a free product [18],
there are but finite number of nonconjugates among $\varp_i$. Q.E.D.

5.4. We move now to the description of homomorphisms of a nonorientable
surface group to a geometrically hyperbolic group $G$.
The key result which will be used in the Classification
Theorem below, is contained in the following lemma.

\proclaim{Lemma} Suppose $G$ is a geometrically hyperbolic group and
$\Sig^g$ is a nonorientable surface, $g>2$.
Then there are, up to an automorphism of $\pi_1(\Sig^g)$ and a conjugation
of $G$, only finitely many homomorphisms $\varp:\pi_1(\Sig^g)\to G$,
such that $\Ker(\varp)$ does not contain a class of a simple closed
essential separating curve.
\endproclaim

\demo{Proof} Let $M$ be as above with $\pi_1(M)=G$.
Fix an orientable double covering $\pi:\bar\Sig\to\Sig$.
Consider a sequence of maps $f_i:\Sig\to M$, which satisfy the condition of the
lemma.
Consider the maps $f_i\circ\pi:\bar \Sig\to M$.
We claim that even though $f_i\circ\pi$ may not be essentially injective,
the proof of 5.3 goes through.
Indeed, to apply the Mumford criterion we need to show that no simple closed
geodesic $\bar C_i$ of $(\bar \Sig,h_i)$ may have a decreasing to zero length.
Observe first that the metrics $h_i$ may be chosen equivariant under the
canonical
involution $\sig$ of $\bar\Sig$, so that $h_i=\pi*\tilde h_i$, where
$\tilde h_i$ is a metric on $\Sig$.
Suppose $length(\bar C_i)\to 0$. The image $C_i=\pi(\bar C_i)$ is a closed
geodesic and $length_{\tilde h_i}(C_i)\to 0$.
Hence there exists a {\bf simple} closed geodesic, which we relabel $C_i$,
such that $length_{\tilde h_i}(C_i)\to 0$.
Since the metric is hyperbolic, it is automatically essential.
Now, either $C_i$ or $C_i^2$ lifts to a simple closed geodesic of $\bar \Sig$
and the argument of [30] shows that $[f_i(C_i^2)]=0$ for $i$ big enough in $G$.
Since $G$ is torsion-free, $C_i\in\Ker f_{i_*}$.
We will show later, that this is impossible.
So the Mumford criterion applies, and we may twist $f_i$ by a diffeomorphism
say $g_i$ of $\bar\Sig$, which by the argument at the end of 6.6 can be taken
$\sig$-equivariant. and get a metric, relabeled $h_i$, which stays in a compact
subset of the Teichmuller space.
Then again as in 5.3, we find a $\sig$-invariant metric $h$ of $\bar \Sig$,
such that the energy of $f_i$ is bounded.
The collar argument of [30] shows that for any simple closed geodesic of
$\bar\Sig$ the minimal length in the homotopy class of its image under $f_i$
stays bounded. This is still true for any simple closed geodesic $C$ in
$\Sig$, since either $C$ or $C^2$ lifts to a simple curve in $\bar\Sig$.
The proof is now completed as in 5.3.

So we need only to show that $f_i$ is essentially injective, i.e.,
$C_i\in\Ker f_{i_*}$ is
 impossible.Relabel $f_i$ by $f$ and $C_i$ by $C$.
Suppose $C\in\Ker f_*$.
By the classification of simple closed curves of [38] we encounter the
following
possibilities (we use the terminology of [38]):
\item{(i)} $[C]=V_1$. Then $V_1^2$ is represented by a simple separating loop
and
$V_1^2\in\Ker f_*$, contrary to the hypothesis of the lemma.

\item{(ii)} $[C]=V_1\ldots V_g$ ($g$ even). This means that after cutting
$\Sig$
along $C$ we get an orientable surface $\hat\Sig$.
Take a separating curve  $D$ in $\hat\Sig$ as in Fig 4.5,left,
then as we have seen, $DCD^{-1}C$ is represented by a simple separating
loop and $DCD^{-1}C\in\Ker f_*$.

\item{(iii)} $[C]=V_1V_2$. Then $V_1^2V_2^2$ is represented by a simple
separating loop and $V_1^2V_2^2\in\Ker f_*$.

\item{(iv)} $[C]=V_1\ldots V_g$ ($g$ odd). Then $[C^2]$ is
represented by a simple
separating loop and $[C^2]\in\Ker f_*$.

\item{(v)} $C$ is separating, which contradicts the hypothesis. Q.E.D.

5.5. We are now ready to state the Classification Theorem for homomorphisms of
a nonorientable surface group to $G$.

\proclaim{Theorem} Let $G$ be a geometrically hyperbolic group, let
$\Sig^g$ be a nonorientable surface and let $\rho:\pi_1\Sig^g)\to G$ be a
homomorphism.
There exists a decomposition
$$
\matrix
\pi_1(\Sig^g)&\longrightarrow&G\hbox{\hskip5truecm}\\
\searrow\del_*&&\uparrow\varp_1*\cdots*\varp_{r+s}*1*\cdots*1\\
&\underbrace{\pi_1(\Sig^{g_1})*\cdots\pi_1(\Sig^{g_r})}_r*
\underbrace{\bbz*\cdots*\bbz}_s*\underbrace{\bbz_2*\cdots\bbz_2}_q&\\
\endmatrix
$$
where $\Sig^{g_1}\cdots\Sig^{g_r}$ are not projective planes, $\del$ is a pinch
and $\varp_1,\ldots,\varp_r$ are essentially injective.
There are, up to automorphisms $\circ f\ \Sig^{g_i}$ and conjugations in $G$,
only finitely
many possibilities for homomorphisms $\varp_i$, $1\le i\le r$.
\endproclaim

\demo{5.6. Example} Let $K\subset S^3$ be an excellent knot.
The knot manifold $S^3\setminus K$ admits a hyperbolic structure of a finite
volume, which does not satisfy 5.1.
Consider the universal cyclic covering $V\to S^3\setminus K$, so that
$\pi_1(V)=G$, where $G$ is the knot group of $K$.
This is a hyperbolic manifold which has just one end, which is a cylinder
$S^1\times\bbr$.
It is easy to check that $V$ satisfies the uniformity and the growth
conditions,
so 5.3 and 5.5 applies for $G$. Recall that G is either f.g. and free or an
infinite
amalgam.

\heading{6. THE DIRICHLET PROBLEM, THE CURRENT NORM AND}\endheading
\heading{THE GENUS ESTIMATES IN}\endheading
\heading{ GEOMETRICALLY HYPERBOLIC GROUPS}\endheading

6.1. In his paper [6], M. Culler proved the following remarkable results for
commutators in free group $F_r$.

\proclaim{6.1.1. Theorem (Culler)} Let $w\in F'_r$ and let $genus(w)$ be the
minimal number of simple commutators, whose product equals $w$. Then
$genus(w^p)$ grows linearly with p.
\endproclaim

\proclaim{6.1.2. Theorem (Culler)} Let $w\in F'_r$ and let $genus (w)=g$.
There are only finite number of equivalence classes of
solutions to the equation
$[x_1,y_1]\cdots[x_g,y_g]=w$ under the action of  $C_w\times A_w$, where the
centrilizer $C_w$ acts by conjugation, and $A_w$ is the subgroup of
$\Aut(F_{2g})$, fixing the word  $r=[x_1,y_1]\cdots [x_g,y_g]$.
\endproclaim

The group $A_w$ is an extension of  the mapping class group  ${\cal M}_g$ [18].

In this section, we will prove theorem 6.2 and below, which extend 6.1.1.
and 6.1.2. to all geometrically hyperbolic groups.

\proclaim{6.2. Theorem} Let $G$ be geometrically hyperbolic and
let $w\in G'$.
Then $genus (w^g)$ grows linearly with $p$.
\endproclaim

\proclaim{6.3. Lemma} Let $M$ be a Riemannian manifold satisfying 5.1.
Then any conjugate class in $\pi_1(M)$ is realized by a closed geodesic.
\endproclaim

\demo{Proof} This is well known (see [9], for example).

\proclaim{6.4. Theorem (Douglas-Rado-Morrey-Lemaire-Jost)}
Let $M$ be a Riemannian manifold satisfying 5.1. and let $\gam$ be an smooth
closed curve in $M$.
There exists a harmonic map of $\Sig^g\setminus D$ in $M$ such that
$\part(\Sig^g\setminus
D)$ is mapped monotonically in $\gam$.
Here $g=genus([\gam])$ and the metric of $\Sig^g\setminus D$ is
chosen arbitrarily.
\endproclaim

\demo{Proof} See [15],[16].

\demo{6.5. Definition} Let $\gam$ be a closed curve in $M$.
Denote $||\gam||_\infty^{-1}$ to be the current norm
$\sup\limits_{\om\in\Om^2(M)}\frac{\int_\gam\om}{||\om||_\infty^{1}}$, where
$||\om||_\infty^1=\sup\limits_M|d\om|$.
\endproclaim

\demo{Proof of the theorem 6.2.} First we choose a closed geodesic $\gam$ with
$[\gam]=\om$ which exists by 6.3.
Let $\Sig^g\setminus D\hookrightarrow M$
be a harmonic map spanning $\gam$, which exists by 6.4.
By the Gauss formula we have $K(\Sig^g\setminus D)\le -k$
in all regular points.
Let $\om\in\Om^2(M)$. Applying the Stocks formula, we get
$$
\int_\gam\om=\int_{\Sig^g\setminus D}d\om\le ||d\om||_\infty\Area(\Sig^g
\setminus D).
$$
The Sharp Thurston Inequality of [25] applies to $\Sig^g\setminus D$ to give
$$
Area(\Sig^g\setminus D)\leq\frac{2\pi}{k}(2g-1)
$$
So
$$
2g-1\ge\frac{k}{2\pi}||\gam||_\infty^{-1}.\tag ****
$$
Applying (****
) to $\gam$, iterated $p$ times, we get
$$
2 genus(\om^p)-1\ge\frac{k\cdot p}{2\pi}||\gam||_\infty^{-1},
$$
as desired.
Observe that yet $p\cdot\gam$, and maybe $\gam$ itself, is not a Jordan
curve, we may slightly perturb it to produce a Jordan curve with arbitrarily
small integral curvature, which is a boundary term in (****) so that the
argument goes
through.

\demo{Proof of the theorem 6.1.1}
Take $M$ to be the universal cyclic covering of the knot manifold of an
excellent
fiber knot, so that $G'$ is finitely generated and hence free [1].
Then $M$ admits a hyperbolic structure as in 5.5. and 6.2. applies directly.

\proclaim{6.6. Theorem} Let $G$ be geometrically hyperbolic and let $w\in G'$
be of
genus $g$.
Then there are only finite number of equivalence classes of solutions
to the equation $[x_1,y_1]\cdots[x_g,y_g]=w$ under the action of
${\cal M}_g$.
\endproclaim

\demo{Proof} Let $M$ be a Riemannian manifold satisfying 5.1. with
$\pi_1(M)=G$.We can always assume $dim(M)\geq 3$.We will provide the argument
under the weak technical condition that $w$ is not a proper power.
Let $\gam$ be a closed geodesic, representing $w$ up to a conjugacy.Perturbing
the metric, make $\gamma$ to be embedded.
Let $f:\Sig^g\setminus D\to M$ be a map, realizing a given solution of
$[x_1,y_1]\cdots[x_g,y_g]=w$ which exists by 1.6.
We can find a conformal structure $C$ on
$\Sig^g\setminus D$ and a harmonic map $\bar f:\Sig^g\setminus D$,
immersion near $\part(\Sig^g\setminus D)$ homotopic to $f$, by [15] and [16].
Moreover, it is easy to find a metric $h$ in $C$, such that the boundary
is geodesic (here we use that $\gam$ is a geodesic in $M$).
So the conformal double $\Sig^{2g}$ is well-defined.
Observe that there exists a conformal involution
$\sig:\Sig^{2g}\to\Sig^{2g}$ with $\part D=\text{Fix}(\sig)$.
We claim that as the solution varies, the conformal structure of $\Sig^{2g}$
can be taken to stay
in a compact set in $M_{12g-6}$.
First observe that
$\Area(f)\le\frac{2\pi(2g-1)}{k}$ by 6.2.
Perturbing the metric of $\Sig$ as in [25] we may assume that
$E^h(\bar f)=\Area(f)\le\frac{2\pi(2g-1)}{k}$.
Adding a strip, we produce a map of $\Sig^{2g}$ to $M$, extending
the map of $\Sig^g\setminus D$ as shown in the Fig 6.6.
An elementary computation shows that the energy of the extended map, say
$\varp$, stays bounded, when $a$ is chosen big enough.

Let $h$ be the unique hyperbolic metric in the conformal class, determined by
this
immersion.
Observe that $\sig$ is necessarily an isometry of $h$.
Consider now the sequence of solutions to
$[x_1,y_1]\cdots[x_g,y_g]=w$ and let $h_i$ be the corresponding
hyperbolic metrics
in $\Sig^{2g}$ and $f_i,\varp_i$ be the corresponding maps.
If $[f_i]$ escapes all compact sets in $M_{12g-6}$, then, as in [30],
we get a closed geodesic $\del_i$ of $(\Sig^{2g},h_i)$ whose length decreases
to
zero, and, again by the collar argument,we can find a curve $\delta'_i$ in the
collar around $\del_i$, such that $\del'_i$ is homotopic to $\del$ and the
image $\bar f_i(\del'_i)$ has the lengh decreasing to zero, and, in particular,
is null-homotopic in $M$.

Consider two cases:
\item{(i)} $\del'_i$ does not intersect $\part D$.
Then we can cut $\Sig^g$ along $\del'_i$ and make the genus reduction
procedure as in 1. which contradicts $g=\text{genus}(w)$.

\item{(ii)} $\del'_i$ intersects $\part D$.
Let $\vare$ be a part of $\del$ which lies in $\Sig^g\setminus \part D$,
with ends in $\part D$, say, $p$ and $q$.
Then we get two different segments $\bar f_i(\vare)$
and a segment $\gam_0$ of $\gam$,
joining $\bar f(p_i)$ and $\bar f(q_i)$, one of which, namely
$\bar f_i(\vare)$, has the length which decreases to zero
and the other lies on the fixed closed geodesic. Hence $length(\gam_0)\to 0$.If
$\e$ and $\gam_0$ are not homotopic (with fixed ends) in $\Sig^g\setminus\part
D$ we repeat the same game: cut $\Sig^g\setminus\part D$ along $\e_i\cup
\gam_0$, glue to $\bar f_i(\e_i\cup\gam_o)$ a disk in $M$ and reduce the genus
of $\gam$ which is impossible. So any $\e_i$ should be homotopic to a part of
$\gam$.The same applies to the mirror of $\Sig^g\setminus D$ and hence
$\del'_i$ is homotopic to an multiple  of $\part D$ which is impossible since
$\bar f_i(\del'_i)$ is null-homotopic in $M$.
So $[f_i]$ stays within the compact set in $M_{12g-6}$.
Twisting by a diffeomorphism
of $\Sig^{2g}$, we can assume that $f_i$ stays in a compact set
$\cale$ in $T_{12g-6}$.
Now, there are only finitely many isotopic classes of conformal involutions
for metrics in $\cale$.

So there are only finite number of isotopy classes of simple curves,
representing $\part D$ after the twist by a diffeomorphism of
$\Sig^{2g}$.
Combining this with the finiteness theorem for $\varp_{i*}:\pi_1(\Sig^{2g})\to
G$,
which follows from the argument of 5.3  complete the proof.

\proclaim{6.7. Corollary} Let $G$ be geometrically hyperbolic, let $w\in G'$
with $genus(w)=g$ and let $A_w$ the stabilizer of $w$ in
$\Aut(G)$.
Then there exists a subgroup of $A_w$ of finite index, which fixes
a subgroup $H\triangleleft G$ of rank $\le 2g$ such that $w\in H'$.
In particular, if $G$ is free and $w$ is a simple commutator,
$w\not=1$, then there exists $N=N(w)$ which that for any endomomorphism
$f$ of $G$, which fixes $w$, $f^N$ fixes one (in fact any)
rank two free group whose commutator contains $w$.
\endproclaim

\demo{Proof} Any endomorphism $f:G\to G$ which fixes $w$, sends solutions of
$[x_1,y_1]\cdots[x_g,y_g]=w$ to solutions.
By 6.6. there are only finite number,
say $N$, of subgroups, generated by $x_i,y_i$, hence a result.

\proclaim{6.8. Corollary} Let $G$ be geometrically hyperbolic, and let
$1\not=u,v\in G'$.
Then there does not exist an endomorphism $\varp$ of $G$ such that
$\varp(u)=u$,
$\varp(v)=uv$.
\endproclaim

\demo{Proof} For $G$ free this is proved in [6], using 6.1.1.
For any $G$, use 6.2. instead of 6.1.1. and complete the proof as in [6].

\centerline{References}

\item{[1]} G.Burde, H.Zieschang, {\sl Knots},de Gruyter, 1985.

\item{[2]} R.G.Burns, C.C.Edmunds, I.H.Farvouqi, {\sl On commutator equations
and stabilazers in free groups}, Canad. Math. Bull, {\bf 19}(1976), 263--267.

\item{[3]} R.G.Burns, C.C.Edmunds, E.Formanek, {\sl The equation
$s^{-1}t^{-1}st=u^{-1}v^{-1}uv$ in free product of groups}, Math. Z., {\bf
153}(1977), 83--88.

\item{[4]} Yu.M.Chmelevsky, {\sl Systems of equations in a free group 1,2},
Soviet Math. Izv., {\bf 35}(1971), 1237--1268; {\bf 36}(1971), 110--179.

\item{[5]} L.P.Commerford, C.C.Edmonds, {\sl Quadratic equations over free
groups and free products},J.of Algebra {\bf 68}(1981),279--297.

\item{[6]} M.Culler, {\sl Using surfaces to solve equations in free groups},
Topology {\bf 20}(1981), 133--145.

\item{[7]} A.L.Edmonds, {\sl Deformations of maps to branched coverings in
dimension two}, Ann.of Math. {\bf 110} (1979), 273--279.

\item{[8]} A.L.Edmonds, I.Bersrtein, {\sl On the construction of branched
covering of low--dimensional manifolds}, Trans. AMS, {\bf 247} (1979), 87--124.

\item{[9]} J.Eells, L.Lemaire, {\sl Another report on harmonic maps},
Bull.Lon.Math.Soc,{\bf 20} (1988), 385--524.

\item{[10]} R.Z.Goldstein, E.C.Turner, {\sl Reduced forms for quadratic words
and spines of 2-manifolds}, Math. Z, {\bf 176}(1981), 195--198.

\item{[11]} R.Grigorchuk, P.Kurchanov, H.Zieschang, {\sl Equivalence of
homomorphisms of surface groups to free groups and some properties of
3-dimensionaal handlebodies}, Proc. Int. Conf. on Algebra, Contemp. Math, {\bf
131}.

\item{[12]} R.Grigorchuk, P.Kurchanov, {\sl Some questions of group theory,
related to geometry}, VINITI Surveys in Modern Math., {\bf 58} (1990),
192--256.

\item{[13]} M.Gromov, {\sl Hyperbolic Groups}, in:Essays in groups theory,
S.M.Gersten, ed.,Springer, 1987, 75--263.

\item{[14]} P.Hill, S.J.Pride, {\sl Commutators, generators and conjugacy
equations in groups}, Arch. Math. {\bf 44}(1985), 1--14.

\item{[15]} J.Jost, {\sl Conformal mappings and the Plateau-Douglas problem in
Riemannian manifolds}, J. Reine Angew. Math, {\bf 359} (1985), 37--54.

\item{[16]} L.Lemaire, {\sl Boundary value problems for harmonic and minimal
maps from surfaces into manifolds}, Ann. Scuola Norm. Super. Pisa, {\bf  9}
(1982), 91--103.

\item{[17]} R.C.Lyndon,{\sl The equation $a^2b^2=c^2$ in free groups}, Michigan
Math. J., {\bf 6} (1959), 155--164.

\item{[18]} R.C.Lyndon, P.E.Schupp, {\sl Combinatorial Group Theory}, Springer,
1977.

\item{[19]} R.C.Lyndon, T.MacDonough, M.Newman, {\sl On powers of products in
groups}, Proc. AMS, {\bf 40},(1973), 419--420.

\item{[20]} R.C.Lyndon, M.J.Wicks, {\sl Commutators in free groups},
Canad.J.Math, {\bf 24}(1981),101--106.

\item{[21]} A.D.Mednych, {\sl Determination of the number of nonequivalent
coverings of a compact Riemann surface}, Sov. Math. Doklady {\bf 239},
(1978),269--271 (Russian).

\item{[22]} D.Piollet, {\sl Solutions d'une equation quadratique dans le groupe
libre}, Diskrete
Math. {\bf 59}(1986), 115--123.

\item{[23]} A.Olshansky, {\sl Diagrams of the surface groups homomorphisms},
Syb. Math. J. {\bf 30}(1989), 150--170.

\item{[24]} A.A.Razborov, {\sl Systems of equations in a free group}, Math.USSR
Izv. {\bf48}(1984), 779--832.

\item{[25]} A.Reznikov, {\sl Harmonic maps, hyperbolic cohomology and higher
Milnor inequalities}, Topology, to appear.

\item{[26]} A.Reznikov, {\sl Yamabe spectra}, preprint (May, 1993).

\item{[27]} G.Rozenberger, {\sl Gleichungen in freien Producten mit Amalgam},
Math. Z.{\bf 173}(1980), 1--12.

\item{[28]} G.Rosenberger, {\sl Producte von Potenzen und Kommutatoren in
freien Gruppen}, J. Algebra,{\bf 53}(1978), 416--422.

\item{[29]} J.Sacks, K.Uhlenbeck, {\sl Minimal immersions of closed Riemann
surfaces}, Trans. AMS, {\bf 271}(1982), 639--652.

\item{[30]} R.Schoen, S.T.Yau, {\sl Existence of incompressible minimal
surfaces and the topology of three  dimensional manifolds with non-negative
scalar curvature}, Ann. Math., {\bf 110}(1979), 127--142.

\item{[31]} P.Schupp, {\sl Quadratic equations in groups cancellation diagrams
on compact surfaces and automorphisms of surface groups}, in:Word Problems 2,
S.I.Adian, W.W.Boone,G.Highman,eds, North-HollandPubl.Comp.,1980, 347--371.

\item{[32]} J.Shapiro, J.Sonn, {\sl Free factor groups of one-relator groups},
Duke Math. J. {\bf 41}(1974), 83--88.

\item{[33]} W.Thurston, {\sl A norm for the homology of 3-manifolds}, Mem.AMS,
{\bf 59} (1986), 98--130.

\item{[34]} W.Thurston, {\sl Three dimensional manifolds, Kleinian groups and
hyperbolic geometry}, Bull. AMS {\bf 6}(1982), 357--381.

\item{[35]} W.Thurston, {\sl Geometry and Topology of Three-Manifolds},
Mimeographed Notes, Princeton, 1979.

\item{[36]} M.Ville, {\sl Harmonic maps, second homology classes of smooth
manifolds and bounded cohomology}, preprint.

\item{[37]} J.A.Wolf, {\sl Spaces of Constant Curvature},McGraw--Hill, 1967.

\item{[38]} H.Zieschang,E.Fogt, H.D.Coldewey,{\sl Surfaces and Planar
Discontinious Groups}, Springer LNM, {\bf 835}, 1980.

\item{[39]} H.Zieschang, {\sl Alternierende Producte in freien Gruppen}, Abh.
Math. Sem. Hamburg, {\bf 27} (1964), 13--31; {\bf 28} (1965), 219--233.

\item{[40]} R. Grigorchuk, P.Kurchanov, {\sl On quadratic equations in free
groups}, Contemp. Math., {\bf 131} (1992), 159--171.

\pbf
\item{} Institute of Mathematics
\item{} The Hebrew University
\item{} Giv'at Ram,91904 Jerusalem
\item{} ISRAEL
\item{}email:simplex@sunrise.huji.ac.il
\bye